\documentclass[12pt]{article}
\usepackage{amsmath, amsthm, amsfonts}
\usepackage{amssymb}
\usepackage{amscd}
\usepackage{verbatim}

\begin{document}
\newcommand{\M}{{\mathcal M}}
\newcommand{\loc}{{\mathrm{loc}}}
\newcommand{\dx}{\,\mathrm{d}x}
\newcommand{\core}{C_0^{\infty}(\Omega)}
\newcommand{\sob}{W^{1,p}(\Omega)}
\newcommand{\sobloc}{W^{1,p}_{\mathrm{loc}}(\Omega)}
\newcommand{\merhav}{{\mathcal D}^{1,p}}
\newcommand{\be}{\begin{equation}}
\newcommand{\ee}{\end{equation}}
\newcommand{\mysection}[1]{\section{#1}\setcounter{equation}{0}}
\newcommand{\bea}{\begin{eqnarray}}
\newcommand{\eea}{\end{eqnarray}}
\newcommand{\bean}{\begin{eqnarray*}}
\newcommand{\eean}{\end{eqnarray*}}
\newcommand{\thkl}{\rule[-.5mm]{.3mm}{3mm}}
\newcommand{\cw}{\stackrel{\rightharpoonup}{\rightharpoonup}}
\newcommand{\id}{\operatorname{id}}
\newcommand{\supp}{\operatorname{supp}}
\newcommand{\wlim}{\mbox{ w-lim }}
\newcommand{\mymu}{{x_N^{-p_*}}}
\newcommand{\R}{{\mathbb R}}
\newcommand{\N}{{\mathbb N}}
\newcommand{\Z}{{\mathbb Z}}
\newcommand{\Q}{{\mathbb Q}}
\newcommand{\abs}[1]{\lvert#1\rvert}
\newtheorem{theorem}{Theorem}[section]
\newtheorem{corollary}[theorem]{Corollary}
\newtheorem{lemma}[theorem]{Lemma}
\newtheorem{definition}[theorem]{Definition}
\newtheorem{remark}[theorem]{Remark}
\newtheorem{proposition}[theorem]{Proposition}
\newtheorem{assertion}[theorem]{Assertion}
\newtheorem{problem}[theorem]{Problem}
\newtheorem{conjecture}[theorem]{Conjecture}
\newtheorem{question}[theorem]{Question}
\newtheorem{example}[theorem]{Example}
\newtheorem{Thm}[theorem]{Theorem}
\newtheorem{Lem}[theorem]{Lemma}
\newtheorem{Pro}[theorem]{Proposition}
\newtheorem{Def}[theorem]{Definition}
\newtheorem{Exa}[theorem]{Example}
\newtheorem{Exs}[theorem]{Examples}
\newtheorem{Rems}[theorem]{Remarks}
\newtheorem{Rem}[theorem]{Remark}

\newtheorem{Cor}[theorem]{Corollary}
\newtheorem{Conj}[theorem]{Conjecture}
\newtheorem{Prob}[theorem]{Problem}
\newtheorem{Ques}[theorem]{Question}
\newcommand{\pf}{\noindent \mbox{{\bf Proof}: }}


\renewcommand{\theequation}{\thesection.\arabic{equation}}
\catcode`@=11 \@addtoreset{equation}{section} \catcode`@=12


\title{On positive solutions of minimal growth
for singular $p$-Laplacian with potential term}
\author{Yehuda Pinchover\\
 {\small Department of Mathematics}\\ {\small
 Technion - Israel Institute of Technology}\\
 {\small Haifa 32000, Israel}\\
{\small pincho@techunix.technion.ac.il}\\\and Kyril Tintarev
\\{\small Department of Mathematics}\\{\small Uppsala University}\\
{\small SE-751 06 Uppsala, Sweden}\\{\small
kyril.tintarev@math.uu.se}}
 \maketitle
\newcommand{\dnorm}[1]{\thkl #1 \thkl\,}

\begin{abstract}
Let $\Omega$ be a domain in $\mathbb{R}^d$, $d\geq 2$, and
$1<p<\infty$. Fix $V\in L_{\mathrm{loc}}^\infty(\Omega)$. Consider
the functional $Q$ and its G\^{a}teaux  derivative $Q^\prime$
given by
$$Q(u)\!:=\!\frac{1}{p}\!\!\int_\Omega\!\!\! (|\nabla u|^p+V|u|^p)\dx,\;\; Q^\prime
(u)\!:=\!-\nabla\cdot(|\nabla u|^{p-2}\nabla u)+V|u|^{p-2}\!u.$$
It is assumed that $Q\geq 0$ on $\core$.  In a previous paper
\cite{ky3} we discussed relations between the absence of weak
coercivity of the functional $Q$ on $\core$ and the existence of a
generalized ground state. In the present paper we study further
relationships between functional-analytic properties of the
functional $Q$ and properties of  positive solutions of the
equation $Q^\prime (u)=0$.
\\[2mm]
\noindent  2000  \! {\em Mathematics  Subject  Classification.}
Primary  \! 35J20; Secondary  35J60, 35J70, 49R50.\\[1mm]
\end{abstract}

\mysection{Introduction}
Properties of positive solutions of quasilinear elliptic
equations, and in particular, of equations with the $p$-Laplacian
term in the principal part defined on general domains have been
extensively studied over the recent decades (see for example,
\cite{AH1,AH2,GS,ky3,TakTin,V} and the references therein).

Fix $p\in(1,\infty)$, a domain $\Omega\subseteq\R^d$ and a
potential $V\in L^\infty_\loc(\Omega)$. So, $\Omega$ is allowed to
be nonsmooth and/or unbounded, and $V$ might blow-up near
$\partial\Omega$ or at infinity.

The $p$-Laplacian equation in $\Omega$ with potential term $V$ is
the equation of the form
\be\label{eq}-\Delta_p(u)+V|u|^{p-2}u=0\quad \mbox{in } \Omega,\ee
where $\Delta_p(u):=\nabla\cdot(|\nabla u|^{p-2}\nabla u)$ is the
$p$-Laplacian. This equation, in the semistrong sense, is a
critical point equation for the functional \be \label{Q}
Q(u)=Q_V(u):=\frac{1}{p}\int_\Omega \left(|\nabla
u|^p+V|u|^p\right)\dx\qquad u\in\core.\ee
So, we consider solutions of \eqref{eq} in the following weak
sense.
\begin{definition}{\em A function
$v\in W^{1,p}_{\mathrm{loc}}(\Omega)$ is a {\em (weak) solution}
of the equation  \be \label{groundstate}
 Q^\prime
(u):=-\Delta_p(u)+V|u|^{p-2}u=0\quad \mbox{in }  \Omega,\ee if for
every $\varphi\in\core$
 \be \label{solution} \int_\Omega (|\nabla v|^{p-2}\nabla
v\cdot\nabla\varphi+V|v|^{p-2}v\varphi)\dx=0. \ee
We say that a real function $v\in C^1_{\mathrm{loc}}(\Omega)$ is a
{\em supersolution} (resp. {\em subsolution})  of the equation
(\ref{groundstate}) if for every nonnegative $\varphi\in\core$
 \be\label{supersolution}
\int_\Omega (|\nabla v|^{p-2}\nabla
v\cdot\nabla\varphi+V|v|^{p-2}v\varphi)\dx\geq 0 \mbox{ (resp.
}\leq 0\mbox{).} \ee
 }\end{definition}
Throughout this paper, unless otherwise stated, we assume that
\be\label{posQ} Q(u)\ge 0 \qquad \forall u\in \core.\ee
The present paper continues the investigation in \cite{ky3} and
also in \cite{aky,TakTin}. These papers deal with global
positivity properties of the functional $Q$ and the set of
positive solutions of the equation $Q'(u)=0$ on a general domain
$\Omega\subset \R^d$. The existence of such global positive
solutions is linked to the positivity of $Q$ by the following
Allegretto-Piepenbrink type theorem.
\begin{theorem}[{\cite[Theorem~2.3]{ky3}}]\label{pos}
Let Q be a functional of the form (\ref{Q}). The following
assertions are equivalent:
\begin{itemize}
 \item[(i)] The functional $Q$ is nonnegative on
$C_0^\infty(\Omega)$.
 \item[(ii)] Equation (\ref{groundstate}) admits a global positive solution.
 \item[(iii)] Equation
(\ref{groundstate}) admits a global positive supersolution.
\end{itemize}
\end{theorem}
It was established in \cite{ky3} that the absence of a weak
coercivity is equivalent to the existence of a (generalized)
ground state. The proof hinged on the representation of $Q$ as an
integral with a nonnegative Lagrangian density due to the
generalized Picone identity \cite{AH1,AH2,DS}. In \cite{aky} a
simplified equivalent (in the sense of two-sided estimates)
expressions for the Lagrangian were derived, and a Liouville-type
theorem was proved. In Theorem~\ref{thmky3} of the present paper,
we prove the equivalence of several weak coercivity properties of
the functional $Q$, as well as their equivalence to the positivity
of the variational $Q$-capacity of closed balls.

For  $p\le d$ it was proved in \cite{ky3} that  the ground state
can be identified as the global positive solution of minimal
growth at infinity of $\Omega$, thus extending {\em criticality
theory} of positive solutions for second-order {\em linear}
elliptic equations (see \cite{P}) to the case of quasilinear
equations of the form \eqref{groundstate}. In
Section~\ref{sec:mingr}, we further study positive solutions of
the equation $Q'(u)=0$ of minimal growth in a neighborhood of
infinity in $\Omega$, and the behavior of positive solutions near
an isolated singularity. Consequently, we extend the above
identification to the case $p>d$ (see
Theorem~\ref{cor_nonremove}).

Finally, we give conditions for a positive solution to be a
positive solution of minimal growth in a neighborhood of infinity
in terms of the infimum of the integral of the corresponding
nonnegative Lagrangian in a neighborhood of infinity. The
variational condition that we present in
Section~\ref{sec:mingr1:2} for the case $p=2$, is necessary and
sufficient, while for $p\neq 2$ we give in
Section~\ref{sec:mingr1} only a stronger sufficient condition.
Weakening this condition with our methods requires the strong
comparison principle, which is generally false for $p\neq 2$.
\mysection{Preliminaries}\label{secprel}
 In this section we  recall {\em
local} properties of solutions of (\ref{groundstate}) that hold in
any smooth subdomain $\Omega'\Subset\Omega$, where $A\Subset B$
means that $\bar{A}$ is compact in $B$.

\vskip 3mm

\noindent {\bf 1. Smoothness and Harnack inequality.}  Weak
solutions of (\ref{groundstate}) admit H\"older continuous first
derivatives, and nonnegative solutions of (\ref{groundstate})
satisfy the Harnack inequality \cite{Serrin1,Serrin2,T}.

\vskip 3mm

\noindent {\bf 2. Harnack convergence principle.} Let
$\{\Omega_{N}\}_{N=1}^{\infty}$ be an {\em exhaustion} of $\Omega$
(i.e., a sequence of smooth, relatively compact domains such that
$\mbox{cl}({\Omega}_{N})\subset \Omega_{N+1}$, and
$\cup_{N=1}^{\infty}\Omega_{N}=\Omega$). Fix a reference point
$x_0\in \Omega_1$. Assume also that  $V,
\{V_{N}\}_{N=1}^{\infty}\subset L^\infty_{\mathrm{loc}}(\Omega)$
satisfy $V_N\to V$ in $L^\infty_{\mathrm{loc}}(\Omega)$. Suppose
that  $\{u_N\}$ is a sequence of positive solutions of the
equations \be\label{Harnackprinc} Q'_N(u_N):
=-\Delta_p(u_N)+V_N|u_N|^{p-2}u_N=0\quad \mbox{in } \Omega_N,\ee
such that $u_N(x_0)=1$.

In light of the Harnack inequality, a priori interior estimates
\cite[Theorem ~1]{T}, the Arzel\`{a}-Ascoli theorem, and a
standard diagonalization argument, there exist $0<\beta<1$, and a
subsequence $\{u_{N_k}\}$ of $\{u_N\}$ that converges in
$C_{\mathrm{loc}}^{1,\beta}(\Omega)$ to a positive solution $u$ of
the equation
$$Q'(u): =-\Delta_p(u)+V|u|^{p-2}u=0\quad \mbox{in } \Omega.$$

\vskip 3mm

\noindent {\bf  3. Principal eigenvalue and eigenfunction.} For
any smooth subdomain $\Omega'\Subset\Omega$ consider the
variational problem \be \label{mu} \lambda_{1,p}(\Omega'):=\inf_{
u\in W_0^{1,p}(\Omega')}\dfrac{\int_{\Omega'}(|\nabla
u|^p+V|u|^p)\dx}{\int_{\Omega'} {|u|^p}\dx}\,.\ee
 It is well-known that for such a subdomain, (\ref{mu}) admits (up to a multiplicative
constant) a unique minimizer $\varphi$ \cite{DKN,GS}. Moreover,
$\varphi$ is a positive solution of the quasilinear eigenvalue
problem
\begin{equation}
  \begin{cases}
Q'(\varphi)=\lambda_{1,p}(\Omega')|\varphi|^{p-2}\varphi    & \text{ in } \Omega', \\
    \varphi=0 & \text{ on } \partial \Omega'.
  \end{cases}
\end{equation}
$\lambda_{1,p}(\Omega')$ and $\varphi$ are called, respectively,
the {\em principal eigenvalue and eigenfunction} of the operator
$Q'$ in $\Omega'$.

It should be noted, however, that minimization statements for
singular elliptic problems on unbounded  or nonsmooth domains
typically do not produce points of minimum, but minimizing
sequences, which locally (up to subsequences) converge to
solutions.

\vskip 3mm

\noindent {\bf 4. Weak and strong maximum principles.}
\begin{theorem}[{\cite{GS} (see also \cite{AH1,AH2})}]\label{thmGS}
Consider a functional $Q$ of the form \eqref{Q} such that
\eqref{posQ} does not necessarily hold in $\Omega$. Let
$\Omega'\Subset\Omega$ be a bounded $C^{1+\alpha}$-subdomain,
where $0<\alpha<1$. So, in particular, $V\in L^\infty(\Omega')$.

\vskip 3mm

The following assertions are equivalent:
\begin{itemize}
 \item[(i)] $Q'$  satisfies the
maximum principle: If $u$ is a solution of the equation
$Q'(u)=f\geq 0$ in $\Omega'$ with some $f\in L^\infty(\Omega')$,
and satisfies $u\geq 0$ on $\partial\Omega'$, then $u$ is
nonnegative in $\Omega'$.

\item[(ii)] $Q'$  satisfies the strong maximum principle: If $u$
is a solution of the equation $Q'(u)=f\gneqq 0$ in $\Omega'$ with
some $f\in L^\infty(\Omega')$, and satisfies $u\geq 0$ on
$\partial\Omega'$, then $u>0$ in $\Omega'$.

\item[(iii)] $\lambda_{1,p}(\Omega')>0$.

\item[(iv)] For some $0\lneqq f\in L^\infty(\Omega')$ there exists a positive strict
supersolution $v$ satisfying  $Q'(v)=f$ in $\Omega'$, and $v=0$ on
$\partial\Omega'$.

\item[(iv$\prime$)] There exists a positive strict supersolution $v$
satisfying   $Q'(v)=f\gneqq 0$ in $\Omega'$, such that $v\in
C^{1+\alpha}(\partial \Omega')$ and $f\in L^\infty(\Omega')$.

\item[(v)] For each nonnegative $f\in C^\alpha(\Omega')\cap
L^\infty(\Omega')$ there exists a unique weak nonnegative solution
of the problem $Q'(u)=f$ in $\Omega'$, and $u=0$ on
$\partial\Omega'$.
 \end{itemize}
 \end{theorem}


\vskip 3mm

\noindent {\bf 5. Weak comparison principle.} We shall need also
the following {\em weak comparison principle} (or WCP for
brevity).
\begin{theorem}[{\cite{GS}}]\label{CP}
Consider a functional $Q$ of the form \eqref{Q} defined on
$\Omega$, such that \eqref{posQ} does not necessarily hold in
$\Omega$. Let $\Omega'\Subset\Omega$ be a bounded subdomain of
class $C^{1,\alpha}$, where $0< \alpha\leq 1$. Assume that
$\lambda_{1,p}(\Omega')>0$ and let $u_i\in W^{1,p}(\Omega')\cap
L^\infty(\Omega')$ satisfying $Q'(u_i)\in L^\infty(\Omega')$,
$u_i|_{\partial \Omega'}\in C^{1+\alpha}(\partial \Omega')$, where
$i=1, 2$. Suppose further that the following inequalities are
satisfied
\begin{equation}
\left\{
\begin{array}{rclc}
 Q'(u_1)\!\!\!&\leq &\!\!\! Q'(u_2) \qquad &\mbox{ in }
\Omega',\\
Q'(u_2)\!\!\!&\geq &\!\!\! 0 \qquad &\mbox{ in }
\Omega',\\
u_{1}\!\!\!&\leq&\!\!\!u_{2} \qquad &\mbox{ on }\partial\Omega',\\
u_{2}\!\!\!&\geq&\!\!\!0 \quad&\mbox{ on } \partial \Omega'.\\
\end{array}\right.
\end{equation}
Then
$$ u_{1}\leq u_{2} \qquad\mbox{ in }\Omega'.$$
\end{theorem}

\vskip 3mm

\noindent {\bf 6. Strong comparison principle.}
\begin{definition}\label{defSCP}{\em we say that the {\em strong
comparison principle} (or SCP for brevity) holds true for the
functional $Q$ if the conditions of Theorem~\ref{CP} imply that
$u_{1}< u_{2}$ in $\Omega'$ unless $u_1=u_2$ in $\Omega'$.}
\end{definition}

\begin{remark}\label{remSCP}{\em
It is well known that the SCP holds true for $p=2$ and for
$p$-harmonic functions. For sufficient conditions for the validity
of the SCP see \cite{AS,CT,DaS,LP,T1} and the references therein.
In \cite{CT} M.~ Cuesta and P.~Tak\'{a}\v{c} present a
counterexample where the WCP holds true but the SCP does not. }
\end{remark}
\mysection{Picone identity and equivalent
Lagrangian}\label{secpicone}
Let $v\in C_{\mathrm{loc}}^1(\Omega)$ be a positive solution
(resp. subsolution) of (\ref{groundstate}). Using the {\em Picone
identity} \cite{AH1,AH2,DS}
 we infer that for every $u\in\core$, $u\ge 0$, we have
 \be\label{QL} Q(u)=\int_\Omega L(u,v)\dx\qquad \left(\mbox{resp. }\; Q(u) \leq \int_\Omega
 L(u,v)\dx\right),
 \ee
where the Lagrangian $L$ is given by \be\label{piconeLag}
L(u,v):=\frac{1}{p}\left[ |\nabla u|^p+(p-1)\frac{u^p}{v^p}|\nabla
v|^p-p\frac{u^{p-1}}{v^{p-1}}\nabla u\cdot|\nabla v|^{p-2}\nabla
v\right].\ee It can be easily verified that $L(u,v)\ge 0$ in
$\Omega$.

Let now $w:=u/v$, where $v$ is a positive solution  of
(\ref{groundstate}) and $u\in\core$, $u\ge 0$. Then \eqref{QL}
implies that
\be \label{QL1} Q(vw)=\frac{1}{p}\int_\Omega\left[|w\nabla
v+v\nabla w |^p-w^p|\nabla v|^p -pw^{p-1}v|\nabla v|^{p-2}\nabla
v\cdot\nabla w\right]\dx. \ee Similarly,  if $v$ is a nonnegative
subsolution of (\ref{groundstate}), then \be \label{QL1sub}
Q(vw)\leq \frac{1}{p}\int_\Omega\left[|w\nabla v+v\nabla w
|^p-w^p|\nabla v|^p -pw^{p-1}v|\nabla v|^{p-2}\nabla v\cdot\nabla
w\right]\dx. \ee
Therefore, a nonnegative functional $Q$ can be represented as the
integral of a nonnegative Lagrangian $L$. But the expression
\eqref{piconeLag} of $L$ contains  an indefinite term, and
consequently the integrand in \eqref{QL1} and \eqref{QL1sub} is
nonnegative  but with nonpositive terms. The next proposition
shows that $Q$ admits a two-sided estimate by a {\em simplified
Lagrangian} containing only nonnegative terms. We call the
functional associated with this simplified Lagrangian the {\em
simplified energy}.

Let $f$ and $g$ be two nonnegative functions. We denote $f\asymp
g$ if there exists a positive constant $C$ such that $C^{-1}g\leq
f \leq Cg$.
\begin{proposition}[{\cite[Lemma~2.2]{aky}}]
\label{prop:superPicone}  Let $v\in C_{\mathrm{loc}}^1(\Omega)$ be
a positive solution of (\ref{groundstate}) and let $w\in
C^1_0(\Omega)$ be a nonnegative function. Then
\bea \label{p<2}  Q(vw)\asymp  \int_\Omega v^2 |\nabla
w|^2\left(w|\nabla v|+v|\nabla w|\right)^{p-2}\dx.\eea
In particular, for $p\ge 2$, we have
\bea \label{p>2}  Q(vw)\asymp  \int_\Omega \left(v^p|\nabla
w|^p+v^2|\nabla v|^{p-2} w^{p-2}|\nabla w|^2\right)\dx. \eea

If $v$ is only a nonnegative subsolution of (\ref{groundstate}),
then
\bea \label{p<2sub}  Q(vw)\leq C  \int_{\Omega\cap\{v> 0\}} v^2
|\nabla w|^2\left(w|\nabla v|+v|\nabla w|\right)^{p-2} \dx.\eea In
particular, for $p\ge 2$, we have
\bea \label{p>2sub}  Q(vw)\leq C  \int_\Omega \left(v^p|\nabla
w|^p+v^2|\nabla v|^{p-2} w^{p-2}|\nabla w|^2\right) \dx. \eea
\end{proposition}
\begin{proof}
Let $1<p<\infty$. The following elementary algebraic vector
inequality holds true \cite{aky}
 \be \label{ineq p 7}
 |a+b|^p-|a|^p-p|a|^{p-2}a\cdot
b\asymp  |b|^2(|a|+|b|)^{p-2} \qquad \forall a,b\in\R^d.\ee Set
now $a:=w|\nabla v|$, $b:=v|\nabla w|$. Then we obtain \eqref{p<2}
and \eqref{p<2sub} by applying
 \eqref{ineq p 7} to \eqref{QL1} and \eqref{QL1sub}, respectively.
\end{proof}
%
\begin{remark}{\em  It is shown in \cite{aky} that for $p>2$ none of the two terms in the
simplified energy \eqref{p>2}  is dominated by the other, so that
\eqref{p>2} cannot be further simplified.}
\end{remark}
\section{Coercivity and ground state}
\begin{definition}\label{defnull}{\em  Let $Q$ be a nonnegative functional on
$\core$. We say that a sequence $\{u_k\}\subset\core$  of
nonnegative functions is a {\em null sequence} of the functional
$Q$ in $\Omega$, if there exists an open set $B\Subset\Omega$ such
that $\int_B|u_k|^p\dx=1$, and \be
\lim_{k\to\infty}Q(u_k)=\lim_{k\to\infty}\int_\Omega (|\nabla
u_k|^p+V|u_k|^p)\dx=0.\ee

We say that a positive function $v\in C^1_{\mathrm{loc}}(\Omega)$
is a {\em ground state} of the functional $Q$ in $\Omega$ if $v$
is an $L^p_{\mathrm{loc}}(\Omega)$ limit of a null sequence of
$Q$.

The functional $Q$ is {\em critical} in $\Omega$ if it admits a
ground state in $\Omega$. If the nonnegative functional $Q$ does
not admit a ground state in $\Omega$, then $Q$ is said to be {\em
subcritical} (or {\em strictly positive}) in $\Omega$.}
\end{definition}
\begin{theorem}[{\cite[Theorem~1.6]{ky3} and \cite{M86,ky2} for the case $p=2$}]
\label{thmgs} Suppose that the functional $Q$ is nonnegative on
$\core$. Then any ground state $v$ is a positive solution of
(\ref{groundstate}). Moreover, $Q$ admits a ground state $v$ if
and only if (\ref{groundstate}) admits a unique positive
supersolution.

In this case, the following Poincar\'e type inequality holds:
There exists a positive continuous function $W$ in $\Omega$, such
that for every $\psi\in C_0^\infty(\Omega)$ satisfying $\int \psi
v \,\mathrm{d}x \neq 0$ there exists a constant $C>0$ such that
the following inequality holds:
 \be\label{Poinc}
  Q(u)+C\left|\int_{\Omega} \psi
 u\,\mathrm{d}x\right|^p\geq C^{-1}\int_{\Omega} W|u|^p\,\mathrm{d}x\qquad
 \forall u\in C_0^\infty(\Omega).\ee
\end{theorem}

The following result extends Theorem~2.7 of \cite{M86}, which was
proved for the linear case, to the case $1<p<\infty$ (cf.
\cite[Theorem~4.2]{ky2}).
\begin{theorem}
\label{loc-unif} Suppose that the functional $Q_V$ is nonnegative
on $\core$. Then $Q_V$ is critical in $\Omega$ if and only if
$Q_V$ admits a null sequence that converges locally uniformly in
$\Omega$.
\end{theorem}
\begin{proof} By Theorem~\ref{thmgs} we only need to show that if
$Q_V$ admits a null sequence, then it admits also a null sequence
that converges locally uniformly in $\Omega$. Let
$\{\Omega_{N}\}_{N=1}^{\infty}$ be an exhaustion of $\Omega$ such
that $x_0\in \Omega_1$. Pick a nonzero nonnegative function $W\in
C_0^\infty(\Omega_1)$. For $t\geq 0$ and $N\geq 1$, consider the
functional $Q_{V-tW}$ on $C_0^\infty(\Omega_N)$. By
Propositions~4.2 and 4.4 of \cite{ky3}, for each $N\geq 1$ there
exists a unique positive number $t_N$ such that the functional
$Q_{V-t_{N}W}$ admits a ground state in $\Omega_N$. Denote by
$v_N$ the corresponding ground state satisfying $v_N(x_0)=1$.

On the other hand, Proposition~4.1 of \cite{ky3} implies that
$t_N\to 0$. Invoking the Harnack convergence principle, it follows
that there exists a subsequence $\{v_{N_k}\}$ of $\{v_{N}\}$ that
converges as $k\to\infty$ locally uniformly to a positive solution
$v$ of the equation $Q_V'(u)=0$ in $\Omega$ such that $v(x_0)=1$.
The uniqueness of a positive solution of the equation $Q_V'(u)=0$
in $\Omega$ satisfying $u(x_0)=1$ (due to the criticality of
$Q_V$) implies that any such subsequence converges to the same
function $v$. Consequently, $v_N\to v$ locally uniformly in
$\Omega$.

Since $v_N\in W^{1,p}_0(\Omega)$, we have
$Q_V(v_N)=p^{-1}t_N\int_\Omega W|v_N|^p\dx$. Consequently,
$$\,Q_V(v_N)=\frac{t_N}{p}\int_\Omega W|v_N|^p\dx\to 0,
\quad \mbox{ and } \int_B |v_N|^p\dx\asymp 1. $$ Therefore,
$\{v_N\}$ is a null sequence of $Q_V$ and  $v$ is the
corresponding ground state.
\end{proof}
A Riemannian manifold $\mathcal{M}$ is said to be {\em
$p$-parabolic} if the equation $$-\Delta_pu=0$$ admits only
trivial positive supersolutions ($p$-superharmonic functions) in
$\mathcal{M}$. In \cite{Tr1,Tr2} Troyanov has established a
relationship between the (variational) $p$-capacity of closed
balls in a Riemannian manifold $\mathcal{M}$ and the
$p$-parabolicity of $\mathcal{M}$. The following is a natural
extension of the definition of $p$-capacity of compact sets.
\begin{definition}[cf. \cite{HKM}] {\em
Suppose that the functional $Q$ is nonnegative on $\core$. Let
$K\Subset \Omega$ be a compact set. The {\em $Q$-capacity} of $K$
in $\Omega$ is defined by
$$\mathrm{Cap}_Q(K,\Omega):=\inf\{Q(u)\mid u\in \core,\; u\geq 1 \;\mbox{ on } K\}.$$
}
\end{definition}
The following theorem extends Theorem~1.6 of \cite{ky3} in the
spirit of \cite[Proposition~3.1]{TakTin}, as well as  Troyanov's
result \cite{Tr1,Tr2} concerning the $p$-capacity of closed balls.
\begin{theorem}
\label{thmky3} Let $\Omega\subseteq\R^d$ be a domain, $V\in
L_\mathrm{loc}^\infty(\Omega)$, and $p\in(1,\infty)$. Suppose that
the functional $Q$ is nonnegative on $\core$. Then the following
statements are equivalent.
\begin{itemize}
\item[(a)] $Q$ does not admit a ground state in $\Omega$.

\item[(b)] There exists a continuous function $W>0$  in $\Omega$
such that
 \be \label{gap_p}
Q(u)\ge \int_\Omega W(x)|u(x)|^p\dx\qquad \forall u\in\core.\ee
\item[(c)]
There exists a continuous function $W>0$  in $\Omega$ such that
\be \label{gap_max} Q(u)\ge \int_\Omega W(x)\left(|\nabla
u(x)|^p+|u(x)|^p\right)\dx \qquad \forall u\in\core. \ee
\item[(d)] There exists an open set $B\Subset\Omega$ and $C_B>0$
such that \be Q(u)\ge C_B\left|\int_B
u(x)\dx\right|^p\label{gap_min}\qquad \forall u\in\core.\ee
\item[(e)] The $Q$-capacity of any closed ball in $\Omega$
is positive.
\end{itemize}
Suppose further that $d>p$. Then $Q$ does not admit a ground state
in $\Omega$ if and only if there exists a continuous function
$W>0$  in $\Omega$ such that \be \label{gap_p*}
Q(u)\ge\left(\int_\Omega W(x)|u(x)|^{p^*}\dx\right)^{p/p^*} \qquad
\forall u\in\core, \ee where $p^*={pd}/{(d-p)}$ is the critical
Sobolev exponent.
\end{theorem}
\begin{proof}  By \cite[Theorem~1.6]{ky3}, (a) $\Leftrightarrow$ (b).
If (b) $\not\Rightarrow$  (c), then there exists a sequence
 $\{u_k\}\subset\core$ with $u_k\geq 0$, and an open set $B\Subset\Omega$ such
that $\int_B|\nabla u_k|^p\dx$=1, while $Q(u_k)\to 0$ and
$\int_B|u_k|^p\dx\to 0$. This is false, since \eqref{QL}, Picone's
formula and Young's inequality imply (see step 1 of the proof of
\cite[Lemma 3.2]{ky3})
$$
\int_B|\nabla u_k|^p\dx\le C Q(u_k)+\int_B|u_k|^p\dx\to 0.
$$
Clearly (c) $\Rightarrow$ (d). If (d) holds and $Q$ admits a
ground state in $\Omega$, then the Poincar\'e type inequality
\eqref{Poinc} implies (b) which is a contradiction.
Statement (e) is immediate from Theorem~\ref{loc-unif}.

If $d>p$, then  \eqref{gap_p*} implies \eqref{gap_min}. On the
other hand,  \eqref{gap_p*} is immediate from \eqref{gap_max} via
partition of unity and the local Sobolev inequality.
\end{proof}
\begin{remark}\label{remc1}{\em
The requirement in Definition~\ref{defnull} that a null sequence
$\{u_k\}$ satisfies  $\{u_k\}\subset \core$ can clearly be
weakened by assuming only that $\{u_k\}\subset W^{1,p}_0(\Omega)$.
Also, the requirement that $\int_B|u_k|^p\dx=1$ can be replaced by
$\int_B|u_k|^p\dx\asymp 1$. Moreover, by Theorem~\ref{thmky3} this
normalization can also be replaced by the requirement that
$\int_Bu_k\dx\asymp 1$.}\end{remark}
\begin{example}\label{ex1}{\em Consider the functional
$Q(u):=\int_{\mathbb{R}^d}|\nabla u|^p\dx$. It follows from
\cite[Theorem~2]{MP} that if $d\leq p$, then $Q$ admits a ground
state $\varphi=\mathrm{constant}$ in $\mathbb{R}^d$. On the other
hand, if $d>p$, then
$$u(x):=\left[1+|x|^{p/(p-1)}\right]^{(p-d)/p}, \qquad v(x):=\mathrm{constant}$$
are two positive supersolutions of the equation $-\Delta_pu=0$ in
$\mathbb{R}^d$. Therefore, $Q$ is strictly positive in
$\mathbb{R}^d$. For further examples see \cite{aky}.
 }\end{example}

\mysection{Solutions of minimal growth at infinity}
\label{sec:mingr} In this section we define and study the
existence of positive solutions of the equation $Q'(u)=0$ of
minimal growth in a neighborhood of infinity in $\Omega$. In
particular, we prove that a positive solution $u$ of the equation
$Q'(u)=0$ in $\Omega$ is a ground state if and only if $u$ is has
minimal growth in any neighborhood of infinity in $\Omega$, a
result which was proved in \cite{ky3} only for $1<p\leq d$.
\begin{definition} {\em
Let $K_0$ be a compact set in $\Omega$.  A positive solution $u$
of the equation $Q'(u)=0$ in $\Omega\setminus K_0$ is said to be a
{\em positive solution of minimal growth in a neighborhood of
infinity in} $\Omega$ (or $u\in\M_{\Omega,K_0}$ for brevity) if
for any compact set $K$ in $\Omega$, with a smooth boundary, such
that $K_0 \Subset \mathrm{int}(K)$, and any positive supersolution
$v\in C((\Omega\setminus K)\cup
\partial K)$ of the equation $Q'(u)=0$ in $\Omega\setminus K$,
the inequality $u\le v$ on $\partial K$ implies that $u\le v$ in
$\Omega\setminus K$.

A (global) positive solution $u$ of the equation $Q'(u)=0$ in
$\Omega$, which has minimal growth in a neighborhood of infinity
in $\Omega$ (i.e. $u\in\M_{\Omega,\emptyset}$) is called a {\em
global minimal solution of the equation $Q'(u)=0$ in $\Omega$}.}
\end{definition}
If $K_0\subset K_1$ are  compact sets in $\Omega$, then clearly
$\M_{\Omega,K_0}\subset \M_{\Omega,K_1}$. On the other hand, the
inverse assertion seems to depend on the SCP. More precisely, we
will prove the following statement after some preparatory lemmas
that are of interest in their own right.
\begin{proposition}
\label{prop_min} Suppose that the functional $Q$ is nonnegative on
$\core$. Assume that the strong comparison principle (SCP) holds
true with respect to $Q$ in any $C^{1,\alpha}$-bounded subdomain
of $\Omega$. Consider two compact sets $K_0,K_1$ in $\Omega$ such
that
$K_0\Subset \mathrm{int}(K_1)$. If $u$
is a positive solution of the equation $Q'(u)=0$ in
$\Omega\setminus K_0$ such that $u\in \M_{\Omega,K_1}$, then $u\in
\M_{\Omega,K_0}$.
\end{proposition}
\begin{definition} {\em Suppose that the functional $Q$ is nonnegative on
$\core$. Let $\{\Omega_{N}\}_{N=1}^{\infty}$ be an exhaustion of
$\Omega$.   Fix $K \Subset \Omega$ with smooth boundary, and let
$u$ be a positive continuous function on $\partial K$. Let $u_N$
be the solution of the following Dirichlet problem
\begin{equation}\label{eq_dirichK}
\left\{
\begin{array}{rcl}
 Q'(u_{N})&=&0 \qquad \mbox{ in }
\Omega_{N}\setminus K,\\
u_{N}&=&u \qquad\mbox{ on }\partial K,\\
u_{N}&=&0 \qquad\mbox{ on }\partial \Omega_{N},
\end{array}\right.
\end{equation}
We denote: \be \label{uK} u^K:=\lim_{N\to\infty} u_N \qquad \mbox{
on } \Omega\setminus K. \ee}
\end{definition}
\begin{lemma}\label{lem_min3}
Suppose that the functional $Q$ is nonnegative on $\core$.  Let
$K$ be a compact set in $\Omega$ with smooth boundary, and let $u$
be a positive continuous function on $\partial K$. Then $u^K$ is
well defined (and in particular, does not depend on the exhaustion
$\{\Omega_{N}\}_{N=1}^{\infty}$). Moreover,  $u^K\in
\M_{\Omega,K}$.
\end{lemma}
\begin{proof}
We note that for $N\geq 1$ we have $\lambda_0(\Omega_{N}\setminus
K)>0$, and therefore the (unique) solvability of
\eqref{eq_dirichK} follows from a standard sub/supersolution
argument and the weak comparison principle (WCP). Using again the
WCP, we see that $\{u_N\}$ is a pointwise nondecreasing sequence.
Moreover, for any $K_0\Subset \mathrm{int}(K)$, and any positive
supersolution $v$ of the equation $Q'(u)=0$ in $\Omega\setminus
K_0$ satisfying $u\leq v$ on $\partial K$, we have by the WCP that
$u_N\leq v$ in $\Omega_{N}\setminus K$. Therefore, the limit $u^K$
exists and $0<u^K\le v$ for any such $v$. Consequently, $u^K$ does
not depend on the exhaustion $\{\Omega_{N}\}_{N=1}^{\infty}$.
Similarly, one checks that $u^K\in \M_{\Omega,K}$.
\end{proof}

\begin{lemma}\label{lem_min} Let $K_0\Subset \Omega$, and let $u$
be a  positive solution of the equation $Q'(u)=0$ in
$\Omega\setminus K_0$. Then $u\in\M_{\Omega,K_0}$ if and only if
for any compact set $K \Subset \Omega$ with smooth boundary, such
that $K_0\Subset \mathrm{int}(K)$, we have $u=u^K$.

Assume further that $K_0$ has a smooth boundary, and that $u$ is
positive and continuous on $(\Omega\setminus K_0)\cup
\partial K_0$. Then $u^{K_0}\leq u$ and equality holds if and only if $u\in\M_{\Omega,K_0}$.
\end{lemma}
\begin{proof}
Let $K\Subset \Omega$ be a set as above, and let $v\in
C((\Omega\setminus K)\cup
\partial K)$ be a positive supersolution of the equation $Q'(u)=0$ in $\Omega\setminus K$,
satisfying the inequality $u\le v$ on $\partial K$. Then by
comparison, $u^K=\lim_{N\to \infty}u_N\leq v$ in $\Omega\setminus
K$. Now if $u=u^K$, it follows that $u\leq v$. Hence
$u\in\M_{\Omega,K_0}$.
\par
On the other hand, if $u\in\M_{\Omega,K_0}$, then by definition
$u\leq u^K$ in $\Omega\setminus K$, and since $u^K\leq u$, we have
$u=u^K$ in $\Omega\setminus K$.

Assume now that $K_0$ has a smooth boundary. Suppose that $u$ is a
positive solution of the equation $Q'(u)=0$ in $\Omega\setminus
K_0$ which is positive and continuous on $(\Omega\setminus
K_0)\cup
\partial K_0$. Then by Lemma~\ref{lem_min3} and its proof we infer that $u^{K_0}\leq u$, and
that equality in this inequality implies that
$u\in\M_{\Omega,K_0}$.

On the other hand, let $u\in\M_{\Omega,K_0}$. Since $u^{K_0}\leq
u$, we need only to prove that $u\leq u^{K_0}$ in $\Omega
\setminus K_0$. In light of the continuity of $u$ and $u^{K_0}$
and the positivity of $u^{K_0}$ we infer that for any
$\varepsilon>0$ there exists a compact smooth set $K_\varepsilon$
satisfying $K_0\Subset \mathrm{int}(K_\varepsilon)$, and
$\mathrm{dist}(\partial K_0,\partial K_\varepsilon)<\varepsilon$,
such that $u\leq (1+\varepsilon)u^{K_0}$ on $\partial
K_\varepsilon$. Since $u\in\M_{\Omega,K_0}$, it follows that
$u\leq (1+\varepsilon)u^{K_0}$ in $\Omega\setminus K_\varepsilon$.
Letting $\varepsilon\to 0$ we obtain $u\leq u^{K_0}$ in $\Omega
\setminus K_0$.
\end{proof}
\begin{remark}\label{rem_min} {\em Let $u$ be a positive solution
as in the first part of Lemma~\ref{lem_min}. Then $u^{K}$ with
$K=K_0$ might be the zero solution (for smooth $K_0$ this happens
if $u|_{K_0}=0$). Therefore, without additional assumptions,  the
set $K$ cannot be replaced in Lemma~\ref{lem_min} by $K_0$.

On the other hand, let $K_0$ be a compact set in $\Omega$ with
smooth boundary, and let $u\in\M_{\Omega,K_0}$ which is positive
and continuous on $(\Omega\setminus K_0)\cup
\partial K_0$.  Then it follows from Lemma~\ref{lem_min} that the comparison
principle for such solutions is also valid on $\Omega\setminus
K_0$ and not only in $\Omega\setminus K$ with  $K_0 \Subset
\mathrm{int}(K)$ as in the definition of $\M_{\Omega,K_0}$.

More precisely,  under the above assumptions on $u$, for any
positive supersolution $v$ of the equation $Q'(u)=0$ in
$\Omega\setminus K_0$ which is positive and continuous on
$(\Omega\setminus K_0)\cup
\partial K_0$ and satisfies $u\leq
v$ on $\partial K_0$, we have $u\leq v$ in $\Omega\setminus K_0$.
 }
\end{remark}

\begin{proof}[Proof of Proposition~\ref{prop_min}] Let $K_0,K_1$ be  compact
sets in $\Omega$ satisfying $K_0\Subset\mathrm{int}(K_1)$, and let
$u$ be a positive solution of the equation $Q'(u)=0$ in
$\Omega\setminus K_0$ such that $u\in \M_{\Omega,K_1}$. Let $K',K$
be smooth compact sets in $\Omega$ satisfying
$$K_0\Subset \mathrm{int}(K')\subset K'\subset K_1\Subset
\mathrm{int}(K)\Subset\Omega.$$ Since $u=u^K$ in $\Omega\setminus
K$ and on $\partial(K\setminus K')$ we have $u\asymp u^{K'}$, it
follows by comparison and exhaustion that  $u\asymp u^{K'}$ in
$\Omega\setminus K'$. Define
$$\varepsilon_{K'}:=\max\{\varepsilon>0\mid \varepsilon u\leq u^{K'} \mbox{ in }
\Omega\setminus K'\}. $$ Clearly, $0<\varepsilon_{K'}\leq 1$.
Suppose that $\varepsilon_{K'}<1$. Since $\varepsilon_{K'}u\lneqq
u^{K'}$ in $\Omega \setminus K'$ and $\varepsilon_{K'}u< u^{K'}$
on $\partial K'$, it follows from the SCP that $\varepsilon_{K'}u<
u^{K'}$ in $\Omega_N \setminus K'$. Therefore, there exists
$\delta>0$ such that $(1+\delta)\varepsilon_0<1$ and
$$(1+\delta)\varepsilon_{K'}u(x)\leq u^{K'}(x) \qquad x\in K\setminus K'.$$
Hence, by comparison and exhaustion argument on $\Omega \setminus
K$, we obtain
$$(1+\delta)\varepsilon_{K'} u(x)=(1+\delta)\varepsilon_{K'} u^K(x)\leq u^{K'}(x)
\qquad x\in \Omega\setminus K.$$ Hence,
$(1+\delta)\varepsilon_{K'} u(x)\leq u^{K'}$ in $\Omega\setminus
K'$ which is a contradiction to the definition of
$\varepsilon_{K'}$. This implies $u\leq u^{K'}$,  and therefore,
$u=u^{K'}$ in $\Omega\setminus K'$. Consequently, by
Lemma~\ref{lem_min}, $u\in M_{\Omega,K_0}$.
\end{proof}
The following two theorems extend (except for the uniqueness
statement of \cite[Theorem~5.4]{ky3}) Theorems~5.4 and 5.5 in
\cite{ky3} which were proved under the assumption $1<p\leq d$.
\begin{theorem}\label{thmmingr1}
Suppose that $1<p<\infty$, and $Q$ is nonnegative on $\core$. Then
for any $x_0\in \Omega$ the equation $Q'(u)=0$ has a positive
solution $u\in\M_{\Omega,\{x_0\}}$.
\end{theorem}
\begin{proof}  Consider an exhaustion
$\{\Omega_{N}\}_{N=1}^{\infty}$ of $\Omega$ such that $x_0\in
\Omega_1$. Let  $\{f_N\}$ be a sequence of nonzero nonnegative
smooth functions such that $f_N$ is compactly supported in
$B(x_0,2/N) \setminus \overline{B(x_0,1/N)}$.

\vspace{3mm}

Fix $N\geq 1$, and denote $A_{N}:=\Omega_N\setminus
\overline{B(x_0,1/N)}$. There exists (with a suitable $c_N>0$) a
unique positive solution of the Dirichlet problem
\begin{equation}
\left\{
\begin{array}{rcl}
 Q'(u_{N})&=&c_Nf_N \qquad \mbox{ in }
A_{N},\\
u_{N}&=&0 \qquad\mbox{ on }\partial A_{N},\\
u_{N}(x_1)&=&1,
\end{array}\right.
\end{equation}
where $x_1\in A_1$ is a fixed reference point. Note that $u_{N}$
is a positive solution of the homogeneous equation in
$\Omega_N\setminus \overline{B(x_0,2/N)}$, and $u_{N}(x_1)=1$. By
the Harnack convergence principle, $\{u_{N}\}$ admits a
subsequence which converges locally uniformly in
$\Omega\setminus\{x_0\}$ to a positive solution $u$ of the
equation $Q'(u)=0$ in $\Omega\setminus\{x_0\}$.
\par
Let $K\Subset \Omega$ be a compact set with a smooth boundary
such that $ x_0\in \mathrm{int}(K)$, and let $v\in
C(\Omega\setminus \mathrm{int}(K))$ be a positive supersolution
of the equation $Q'(u)=0$ in $\Omega\setminus K$
 such that the inequality $u\le v$
holds on $\partial K$.
\par
For  $N\geq N_K$ we have that $\supp f_N\subset B(x_0,2/N) \Subset
K$. Fix $\delta
>1$ and $x\in \Omega\setminus K$.  Applying the WCP
(Theorem~\ref{CP}) in $\Omega_N\setminus K$, we obtain that
$u_N(x)\leq (1+\delta) v(x)$ for all $N$ sufficiently large. By
letting $N\to \infty$, and then $\delta \to 0$, we obtain that
$u\leq v$ in $\Omega\setminus K$. Hence,
$u\in\M_{\Omega,\{x_0\}}$. \vspace{3mm}
\end{proof}
\begin{theorem}\label{thmmingr13}
The functional $Q_V$ is strictly positive in $\Omega$ if and only
if the equation $Q_V'(u)=0$ does not admit a global minimal
solution in $\Omega$. In particular, $u$ is ground state of the
equation $Q_V'(u)=0$ in $\Omega$ if and only if $u$ is a global
minimal solution of this equation.
\end{theorem}
\begin{proof}
{\em 1. Necessity.} Assume that there exists a global minimal
solution $u>0$ of the equation $Q_V'(u)=0$ in $\Omega$, and
suppose that $Q_V$ is strictly positive. It follows from
Theorem~\ref{thmky3} and Theorem~\ref{pos} (see also
\cite[Proposition~4.4]{ky3}) that there exists a nonzero
nonnegative function $V_1\in \core$ with $\supp V_1\subset
B(x_0,\delta)$ for some $\delta>0$, such that $Q_{V-V_1}$ is
strictly positive in $\Omega$. Therefore, due to
Theorem~\ref{pos}, there exists a positive solution $v$ of the
equation $Q'_{V-V_1}(u)=0$ in $\Omega$.
\par
In particular, $v$ is a positive supersolution of the equation
$Q'_V(u)=0$ in $\Omega$ which is not a solution. On the other
hand, $u$ is a positive solution of  the equation $Q'_V(u)=0$ in
$\Omega$ which has minimal growth in a neighborhood of infinity in
$\Omega$. Therefore, there exists $\varepsilon>0$ such that
$\varepsilon u\leq v$ in $\Omega$. Define
$$\varepsilon_0:=\max\{\varepsilon>0\mid \varepsilon u\leq v \mbox{ in } \Omega\}. $$
Clearly $\varepsilon_0 u\lneqq v$ in $\Omega$. Consequently, there
exist $\delta_1,\delta_2>0$ and $x_1\in \Omega$ such that
$$(1+\delta_1)\varepsilon_0 u(x)\leq v(x) \qquad x\in B(x_1,\delta_2).$$
Hence, by the definition of minimal growth, we have
$$(1+\delta_1)\varepsilon_0 u(x)\leq v(x) \qquad x\in \Omega\setminus B(x_1,\delta_2),$$
and thus $(1+\delta_1)\varepsilon_0 u\!\leq\!v$ in $\Omega$, which
is a contradiction to the definition of $\varepsilon_0$.

\vspace{3mm}

{\em 2. Sufficiency.} Fix $x_0\in \Omega_1$. Assume that $Q$ is
not strictly positive. Then $Q$ admits a (unique) ground state $u$
in $\Omega$ satisfying $u(x_1)=1$, where $x_1\in
\Omega\setminus\Omega_1$ is another fixed reference point. It
suffices to prove that $u$ is a global minimal solution of the
equation $Q'(u)=0$ in $\Omega$.
\par
Fix $n\in \mathbb{N}$, and let $f_n\in C_0^\infty(B(x_0,1/n))$ be
a nonzero nonnegative function. For $N\geq 1$, let $v_{N,n}$ be
the unique positive solution of the Dirichlet problem
\begin{equation}
\left\{
\begin{array}{rcl}
 Q'(v_{N,n})&=&f_n \qquad \mbox{ in }
\Omega_{N},\\
v_{N,n}&=&0 \qquad\mbox{ on }\partial\Omega_{N}.
\end{array}\right.
\end{equation}
By the WCP, $\{v_{N,n}\}_{N\geq 1}$ is a nondecreasing sequence.
Recall from \cite[Theorem 1.6]{ky3} that if $Q$ admits a ground
state $u$ in $\Omega$, then any positive supersolution for the
equation $Q'(u)=0$ in $\Omega$ equals to $cu$ with some $c>0$. On
the other hand, if $\{v_{N,n}(x_1)\}$ is bounded, then $v_{N,n}\to
v_n$, where $v_n$ satisfies the equation $Q'(v_n)=f_n\gvertneqq 0$
in $\Omega$. Therefore, $v_{N,n}(x_1)\to\infty$.

Consider now the sequence of functions
$u_{N,n}(x):=v_{N,n}(x)/v_{N,n}(x_1)$, $N\geq 1$. Then  $u_{N,n}$
solves the Dirichlet problem
\begin{equation}
\left\{
\begin{array}{rcl}
 Q'(u_{N,n})&=&\displaystyle\frac{f_n(x)}{v_{N,n}(x_1)^{p-1}} \qquad \mbox{ in }
\Omega_{N},\\[4mm]
u_{N,n}&=&0 \qquad\qquad\quad\mbox{ on }\partial\Omega_{N},\\[2mm]
 u_{N,n}(x_1)&=&1.
\end{array}\right.
\end{equation}
 By the Harnack convergence principle, we
may extract a subsequence of $\{u_{N,n}\}$ that converges as
$N\to\infty$ to a positive solution $u_n$ of the equation
$Q'(u)=0$ in $\Omega$. By the uniqueness of the ground state, we
have $u_n=u$.

Let $K\Subset \Omega$ be a compact set with a smooth boundary, and
let $v\in C(\Omega\setminus \mathrm{int}(K))$ be a positive
supersolution of the equation $Q'(u)=0$ in $\Omega\setminus K$
such that the inequality $u\le v$ holds on $\partial K$. Without
loss of generality, we may assume that $ x_0\in \mathrm{int}(K)$.
\par
Let $n\in \mathbb{N}$ be sufficiently large number  such that
$\supp f_n \Subset K$. By comparison it follows (as in the first
part of the proof) that $u=u_n\leq v$ in $\Omega\setminus K$.
Since $K\Subset \Omega$ is an arbitrary  smooth compact set, it
follows that the ground state $u$ is a global minimal solution of
the equation $Q'(u)=0$ in $\Omega$.
\end{proof}
Consider a positive solution $u$ of the equation $Q'(u)=0$ in a
punctured neighborhood of $x_0$ which has a nonremovable
singularity at $x_0\in \mathbb{R}^d$. Without loss of generality,
we may assume that $x_0=0$. If  $1<p\leq d$, then the behavior of
$u$ near an isolated singularity is well understood.  Indeed, due
to a result of L.~V\'{e}ron  (see \cite[Lemma~5.1]{ky3}), we have
that
 \begin{equation}\label{nonremovasymp}
  u(x)\sim\begin{cases}
    \abs{x}^{\alpha(d,p)} & p<d, \\
     -\log \abs{x} & p=d,
  \end{cases} \qquad \mbox{ as } x\to 0,
\end{equation}
where $\alpha(d,p):=(p-d)/(p-1)$, and $f\sim g$ means that $$
\lim_{x\to 0}\frac{f(x)}{g(x)}= C$$ for some positive constant
$C$. In particular,  $\lim_{x\to 0} u(x)=\infty$.

Assume now that $p>d$. A general question is whether in this case,
any positive solution of the equation $Q'(u)=0$ in a punctured
ball centered at $x_0$ can be continuously extended at $x_0$ (see
\cite{M} for partial results).

We answer below this question under the assumption that  $u\asymp
1$ near the isolated singular point.

\begin{lemma}\label{lemveron}
Assume that $p> d$, and let  $v$ be a positive solution of the
equation $Q'(u)=0$ in a punctured neighborhood of $x_0$ satisfying
$v\asymp 1$ near $x_0$.
 Then $v$ can be continuously extended at $x_0$.
\end{lemma}
\begin{proof}
We use V\'{e}ron's  method (see \cite[Lemma~5.1]{ky3}). Without
loss of generality, we may assume that $x_0=0$. For $0<r\leq r_0$
denote
$$ m(r):=\min_{|x|=r}v(x),\quad M(r):=\max_{|x|=r}v(x), \quad M:=\limsup_{r\to 0} M(r).$$
By our assumption $0< M<\infty$.
 Let $\{x_n\}$ be a sequence such that $x_n\to 0$, and
$M =\lim_{n\to\infty}v(x_n)$. Define $v_n(x):=v(\abs{x_n}x)$.

Since $u\asymp 1$ near $x_0$, it follows that there exists $C>0$
such that in an arbitrarily large punctured ball
$$C^{-1}\leq v_n(x)\leq C$$  for all $n$ large enough.
Moreover, in such a ball, $v_n$ is a positive solution of the
quasilinear elliptic equation
 \be\label{eqdilation}-\Delta _p
v_n(x) +\abs{x_n}^pV(\abs{x_n}x)v_n^{p-1}(x)=0.\ee


 Since $\{v_n\}$ is locally bounded and bounded away
from zero in any punctured ball, the Harnack convergence principle
implies that there is a subsequence of $\{v_n\}$ that converges
locally uniformly in $\mathbb{R}^d\setminus \{0\}$ to a positive
bounded solution $U$ of the limiting equation $-\Delta _p U=0$ in
the punctured space. Recall that by Example~\ref{ex1}, if $p\geq
d$, then the constant function is a ground state of the
$p$-Laplacian \cite{aky}, and in particular, it is the unique
positive $p$-(super)harmonic function in $\mathbb{R}^d$.
Therefore, Theorem~\ref{thmky3} implies that the origin has zero
$p$-capacity. By Theorem~7.36 in \cite{HKM}, the singularity of
$U$ at the origin is removable. Hence, $U$ is an entire positive
$p$-harmonic function in $\mathbb{R}^d$, and consequently,
$U=\mathrm{constant}=\alpha$.

 This implies that
\begin{equation}\label{eqvmu}
\lim_{n\to\infty} \| {v(x)-\alpha}\|_{L^\infty(\abs x=\abs{
x_n})}=0.\end{equation} In other words, $v$ approximates $\alpha$
uniformly over concentric spheres whose radii converge to $0$.

We need to  prove that $v$ approximates $\alpha$ uniformly over
the concentric annuli
 $A_n:=\{\abs{ x_n}\leq \abs{x}\leq \abs{
x_{n+1}}\}$. Let $\alpha_-(x):=\alpha-\delta \abs x^{a}$ and
$\alpha_+(x):=\alpha+\delta \abs x^{a}$ (for some $a>0$ and
$\delta>0$ sufficiently small). It turns out that $\alpha_-$
(resp., $\alpha_+$) is a radial positive subsolution (resp.,
supersolution) of the equation $Q'(u)=0$ near the origin, and
therefore using the comparison principle in the annulus $A_n$, $n\in\N$, and
\eqref{eqvmu},  it follows that
$$\lim_{r\to 0} \| {v(x)-\alpha}\|_{L^\infty(\abs x=
r)}=0.$$
\end{proof}
Let $u\in\M_{\Omega,\{0\}}$ (without loss of generality, we assume
that $x_0=0\in \Omega$). If $u$ has a removable singularity at
$0$, then by definition, $u$ is a global minimal solution. Let us
show that the converse is also true.

Suppose that $u$ has a nonremovable singularity at $0$, and set
$$m\!: =\!\liminf_{r\to 0} m(r)\!=\! \liminf_{r\to 0}\min_{|x|=r} u(x),\quad
M\!: =\!\limsup_{r\to 0} M(r)\!=\! \limsup_{r\to 0}\max_{|x|=r}
u(x).$$

By Harnack inequality, for any positive solution $v$ of the
equation $Q'(u)=0$ in a punctured neighborhood of $0$  there
exists $\tilde{C}>0$ such that
\be\label{eqHarnack}\tilde{C}^{-1}M(r)\leq m(r)\leq M(r) \qquad
0<r\leq r_0.\ee

If $m=0$, then by comparing $u$ with any positive global
(super)solution and using Harnack inequality \eqref{eqHarnack}, we
infer that $u$ must be identically zero, which is a contradiction.

On the other hand, if $M=\infty$, then the equation $Q'(u)=0$ in
$\Omega$ does not admit a global minimal solution due to the
Harnack inequality \eqref{eqHarnack} and the weak comparison
principle. Hence, Theorem~\ref{thmmingr13} implies that $Q$ is
strictly positive.

Assume now that $0<m\leq M<\infty$. Then by Lemma~\ref{lemveron},
$u$ can be continuously extended at the origin. If the equation
$Q'(u)=0$ in $\Omega$ admits a global minimal solution $v$, then
by comparison, $u=cv$, where $c$ is a positive constant, and thus
$u$ has a removable singularity at $0$, a contradiction.

Thus, we proved the following result which
extends the second part of \cite[Theorem~5.4]{ky3}, where the case
$1<p\leq d$ is considered.
\begin{theorem}\label{cor_nonremove} Let $x_0\in \Omega$, and let $u\in\M_{\Omega,\{x_0\}}$.
Then $Q$ is strictly positive in $\Omega$ if and only if $u$ has a
nonremovable singularity at $x_0$.
\end{theorem}
\mysection{Variational principle for solutions of minimal growth
in the linear case} \label{sec:mingr1:2} Throughout this section
we restrict our consideration to the linear case ($p=2$). In fact,
as in \cite{ky2},  we can actually consider in the linear case the
following somewhat more general functional than $Q_V$ of the form
\eqref{Q}.

We assume that $A:\Omega \rightarrow \mathbb{R}^{d^2}$ is a
measurable symmetric matrix valued function such that for every
compact set $K\Subset \Omega$ there exists $\mu_K>1$ so that \be
\label{stell} \mu_K^{-1}I_d\le A(x)\le \mu_K I_d \qquad \forall
x\in K,
 \ee
where $I_d$ is the $d$-dimensional identity matrix, and the matrix
inequality $A\leq B$ means that $B-A$ is a nonnegative matrix on
$\mathbb{R}^d$. Let $V\in
L^{q}_{\mathrm{loc}}(\Omega;\mathbb{R})$, where $q>{d}/{2}$. We
consider the quadratic form \be \label{assume}
\mathbf{a}_{A,V}[u]:=\frac{1}{2}\int_\Omega\left(A\nabla u\cdot
\nabla u+V|u|^2\right)\mathrm{d}x   \ee on $\core$
 associated with the Schr\"odinger
equation \be \label{divform}
Pu:=(-\nabla\cdot(A\nabla)+V)u=0\qquad \mbox{ in } \Omega. \ee
 We say that $\mathbf{a}_{A,V}$ is {\em nonnegative} on $C_0^\infty(\Omega)$,
 if $\mathbf{a}_{A,V}[u]\geq 0$ for all $u\in \core$.

Let $v$ be a positive solution of the equation $Pu=0$ in $\Omega$.
Then by \cite[Lemma~2.4]{ky2} we have the following analog of
\eqref{QL}. For any nonnegative $w\in \core$ we have
 \begin{equation}
\label{e1:2} \mathbf{a}_{A,V}[w]= \frac{1}{2}\int_\Omega v^2
A\nabla (w/v)\cdot\nabla (w/v) \,\mathrm{d}x.
 \end{equation}
Moreover, it follows from \cite{ky2,ky3} that all the results
mentioned in the present paper concerning
criticality/subcriticlity of the functional $Q$ are also valid for
the form $\mathbf{a}_{A,V}$.

The aim of this section is to characterize positive solutions of
minimal growth in a neighborhood of infinity in $\Omega$ in terms
of a modified null sequence of the form $\mathbf{a}_{A,V}$.

\begin{theorem}\label{thmmin_null:2}
Suppose that $\mathbf{a}_{A,V}$ is nonnegative on $\core$. Let
$\Omega_1\Subset \Omega$ be an open set, and let $u\in
C(\Omega\setminus \Omega_1)$ be a positive solution of the
equation $Pu=0$ in $\Omega\setminus \overline{\Omega_1}$.
\par
Then $u\in\M_{\Omega,\overline{\Omega_1}}$ if and only if for
every smooth open set $\Omega_2$ satisfying $\Omega_1 \Subset
\Omega_2\Subset \Omega$, and an open set $B\Subset(\Omega\setminus
\overline{\Omega_2})$ there exists a sequence
$\{u_k\}\subset\core$,  $u_k\ge 0$, such that for all $k\in\N$,
$\int_B |u_k|^2\dx=1$, and \be
\lim_{k\to\infty}\int_{\Omega\setminus \overline{\Omega}_2} u^2
A\nabla (u_k/u)\cdot\nabla (u_k/u)\dx =0.\ee
\end{theorem}
\begin{proof}
{\it 1. Sufficiency.} Let $u\in C(\Omega\setminus \Omega_1)$ be
a positive solution of the equation $Pu=0$ in $\Omega\setminus
\overline{\Omega}_1$.  Let $\Omega_2$ be an open set with smooth
boundary such that $\Omega_1 \Subset \Omega_2\Subset \Omega$, and
let $B$ be an open set so that $B\Subset(\Omega\setminus
\overline{\Omega_2})$.
\par
Suppose that $\{u_k\}\subset\core$ is a sequence of nonnegative
functions such that $\int_B|u_k|^2\dx=1$ for all $k\in\N$, and
\be\label{eq4.6:2} \lim_{k\to\infty}\int_{\Omega\setminus
\overline{\Omega_{2}}}u^2 A\nabla (u_k/u)\cdot\nabla (u_k/u)\dx
=0.\ee
It follows (cf. \cite[Lemma~2.5]{ky2}) that $u_k\to cu$ in
$W^{1,2}_{\mathrm{loc}} (\Omega\setminus\overline{\Omega_{2}})$,
where $c>0$.

Let $E:W^{2,2}_{\mathrm{loc}}(\Omega \setminus \Omega_{2}) \to
W^{2,2}_{\mathrm{loc}}(\Omega)$ be a bounded linear extension
operator (cf. \cite[Section~5.4]{E}, and in particular, Remark (i)
therein). Note that the operator $E$ in \cite{E} extends a given
function to a bounded open set outside a smooth compact boundary,
so the construction from \cite{E} applies also to our situation,
combined with a straightforward use of partition of unity. For
completeness, we outline the construction of $E$ below.

In suitable neighborhoods $U_1,\dots,U_m$ covering
$\partial\Omega_{2}$ there exist local diffeomorphisms
$\{\psi_j\}_{1\leq j\leq m}$, such that for each $1\leq j\leq m$
the diffeomorphism $\psi_j$ maps the open set
$U_j\cap(\Omega\setminus\overline{\Omega_{2}})$ into a portion of
the half-space $\mathbb{R}^d_+=\{(x',x_d)\in \mathbb{R}^d \mid
x_d>0\}$, and $\psi_j(U_j\cap(\partial \Omega_{2}))\subset
\{(x',0)\in \mathbb{R}^d\}$. We may assume that $w\in
C^\infty(\Omega \setminus \Omega_{2})$ and define in
$U_1,\dots,U_m$, up to the diffeomorphism
$$w_j(x',x_d):=
  \begin{cases}
    w(x',x_d) & \text{ if } x_d\geq 0, \\
   4w(x',-x_d/2)-3w(x',-x_d)  & \text{ if } x_d<0.
  \end{cases}$$
  Let $U_0$
and $U_\infty$ be open sets, $U_0\Subset\Omega_{2}$ and
$\overline{U_\infty} \subset\Omega\setminus\overline{\Omega_{2}}$,
such that together with $U_1,\dots,U_m$ they form an open covering
of $\Omega$. Set $w_0(x):=1$ for  $x\in U_0$, and
$w_\infty(x):=w(x)$ if $x\in U_\infty$. Let
$\{\chi_j\}_{j=0,\ldots,m,\infty}$ be a partition of unity
subordinated to $U_0,\dots,U_m,U_\infty$, and define the extension
operator as
$$
(Ew)(x):=w_\infty(x)\chi_\infty(x)+\sum_{j=0}^mw_j(x)\chi_j(x).
$$
Let $\tilde{u}:=Eu$, and $\tilde{v}_k:=E(u_k/u)$, and
$\tilde{u}_k:=\tilde{u}\tilde{v}_k$.
 %
 %
Note that, since $u>0$ and $\partial\Omega_2$ is compact, one can always choose the neighborhoods
$U_1,\dots,U_m$ sufficiently small so that $\tilde{u}>0$. Clearly,
$Ew|_{\Omega\setminus{\overline \Omega_2}}=w$. Let
$f:=P\tilde{u}$, and define $W:=f/\tilde{u}$. Notice that $W$ has
a compact support in $\Omega$. Moreover,  by elliptic regularity,
$W\in L^q(\Omega)$. It follows that $\tilde{u}$ is a positive
solution of the equation $(P-W)u=0$ in $\Omega$. Moreover, by the
continuity of $E$ and the continuity of  the form $\mathbf{a}$ due
\eqref{e1:2}, it follows from \eqref{eq4.6:2} that \be
\lim_{k\to\infty} \mathbf{a}_{A,V-W}[\tilde{u}_k]=
 %
 %
\lim_{k\to\infty}\int_{\Omega} \tilde{u}^2 A\nabla
(\tilde{u}_k/\tilde{u})\cdot\nabla (\tilde{u}_k/\tilde{u})
\,\mathrm{d}x=0.\ee On the other hand, $\int_B|\tilde{u}_k|^2\dx
=\int_B|u_k|^2\dx= 1$. Therefore,  Corollary ~1.6 in \cite{ky2}
implies that $\tilde{u}_k\to c\tilde{u}$, and
$\tilde{u}\in\M_{\Omega,\emptyset}$ (Corollary~1.6 of \cite{ky2}
is analogous to Theorem~\ref{thmmingr13}, but note that the
terminology in \cite{ky2} is different from the terminology of the
present paper). Hence $u\in \M_{\Omega,\overline{\Omega_2}}\,$.
Since $\Omega_2$ is an arbitrary smooth open set satisfying
$\Omega_1 \Subset \Omega_2\Subset \Omega$, we have $u\in
\M_{\Omega,\overline{\Omega_1}}\,$.

\vskip 1.2mm

{\it 2. Necessity.} Suppose that  $u\in
C(\Omega\setminus\Omega_1)\cap \M_{\Omega,\overline{\Omega_1}}$.
Let $\Omega_2$ be any open set with smooth boundary such that
$\Omega_1 \Subset \Omega_2\Subset \Omega$.

Let $\tilde{u}$ be a positive function in
$W^{1,2}_{\mathrm{loc}}(\Omega)$ such that
$\tilde{u}|_{\Omega\setminus\Omega_2}=u$. Let $f:=P\tilde{u}$, and
define $W:=f/\tilde{u}$. Recall that $W$ has a compact support in
$\Omega$, and that the SCP holds in the linear case. Therefore,
Proposition~\ref{prop_min} implies that $\tilde{u}$ is a global
minimal solution of the equation $(P-W)u=0$ in $\Omega$.
Consequently, it follows from \cite[Corollary~1.6]{ky2} that
$\tilde{u}$ is a ground state of the equation $(P-W)u=0$ in
$\Omega$. Let $\{u_k\}\subset \core$ be a null sequence for the
form $\mathbf{a}_{A,V-W}$. So, for some open set
$B\Subset(\Omega\setminus \overline{\Omega_2})$ we have
$\int_B|u_k|^2\dx= 1$, and
$$\lim_{k\to\infty}\int_{\Omega} \tilde{u}^2 A\nabla
({u}_k/\tilde{u})\cdot\nabla ({u}_k/\tilde{u})
\,\mathrm{d}x=\lim_{k\to\infty} \mathbf{a}_{A,V-W}[u_k] =0.$$
Thus, \be \lim_{k\to\infty}\int_{\Omega\setminus
\overline{\Omega}_2} {u}^2 A\nabla ({u}_k/u) \cdot\nabla ({u}_k/u)
\,\mathrm{d}x =0.\ee
\end{proof}

Finally, we prove a sub-supersolution comparison principle for our
singular elliptic equation. This general Phragm\'{e}n-Lindel\"{o}f
type principle, which seems to be new even in the linear case,
holds in unbounded or nonsmooth domains, and for irregular
potential $V$, provided the subsolution satisfies a certain decay
property of variational type (cf. \cite{Agmon84,LLM,MMP,PS}).
\begin{theorem}[Comparison Principle]\label{thm_comparison:2}
Let $P$ be a nonnegative Schr\-\"odinger operator of the form
\eqref{divform}. Fix smooth open sets $\Omega_1\Subset
\Omega_2\Subset \Omega$. Let $u,v\in
W^{1,2}_{\mathrm{loc}}(\Omega\setminus \Omega_1)\cap
C(\Omega\setminus \Omega_1)$ be, respectively, a positive
subsolution and a supersolution of the equation $Pw=0$ in
$\Omega\setminus \overline{\Omega_1}$ such that $u\leq v$ on
$\partial \Omega_2$.
\par
Assume further that $Pu\in L^q _{\mathrm{loc}}(\Omega\setminus
\Omega_1)$, where $q>{d}/{2}$, and that there exist an open set
$B\Subset(\Omega\setminus \overline{\Omega_2})$ and a sequence
$\{u_k\}\subset\core$, $u_k\ge 0$, such that \be\label{condsubsol}
\int_B |u_k|^2\dx=1 \quad \forall k\geq 1, \mbox{ and }\quad
\lim_{k\to\infty}\int_{\Omega\setminus \overline{\Omega}_1} {u}^2
A\nabla ({u}_k/u)\cdot\nabla ({u}_k/u) \,\mathrm{d}x =0.\ee Then
$u\leq v$ on $\Omega\setminus \Omega_2$.
\end{theorem}
\begin{proof} Let $\Omega_1\Subset
\Omega_2\Subset \Omega$, and let $0<\tilde{u}\in
W^{1,2}_{\mathrm{loc}}(\Omega)$  be the extension  of $u$ provided
in the proof of Theorem~\ref{thmmin_null:2}, and define
analogously $f:=P\tilde{u}$ and $W:=f/\tilde{u}$. Clearly,
$W\leq 0$ in $\Omega\setminus \Omega_2$.
Theorem~\ref{thmmin_null:2} and Proposition~\ref{prop_min} imply
that $\tilde{u}$ is a positive solution of the equation $(P-W)u=0$
in $\Omega$ which is a global minimal solution. On the other hand,
$v$ is a positive supersolution of the equation $(P-W)u=0$ in
$\Omega\setminus \Omega_2$, therefore $u\leq v$ on
$\Omega\setminus \Omega_2$.
\end{proof}
\begin{remark}\label{remsub:2}
{\em In Theorem~\ref{thm_comparison:2} we have assumed that the
subsolution $u$ is strictly positive. It would be useful to prove
the above comparison principle under the assumption that $u\geq 0$
(cf. \cite{LLM}).}
\end{remark}
\begin{remark}\label{remsub:3}
{\em Let $K$  be a compact set in $\Omega$, and $\phi$ be a
positive solution of minimal growth in a neighborhood of infinity
in $\Omega$ of the equation
$\tilde{P}u:=-\nabla\cdot(\tilde{A}\nabla u)+\tilde{V}u=0$ in
$\Omega$ for some $\tilde{A}$ satisfying \eqref{stell}, and
$\tilde{V}\in L^{q}_{\mathrm{loc}}(\Omega;\mathbb{R})$, with
$q>{d}/{2}$. If $u\in W^{1,2}_{\mathrm{loc}}(\Omega\setminus
\Omega_1)\cap C(\Omega\setminus \Omega_1)$ is a positive
subsolution of the equation $Pw=0$ in $\Omega\setminus
\overline{\Omega_1}$ such that
$$u^2(x)A(x)\leq \phi^2(x)\tilde{A}(x)\qquad  \mbox{ in } \Omega\setminus K,$$
then Condition \eqref{condsubsol} is satisfied (cf.
\cite{PLiouv}).}
\end{remark}


\mysection{Variational principle for solutions of minimal growth
for the quasilinear case} \label{sec:mingr1} In this section we
extend the results of the previous section to the case
$1<p<\infty$. Since the SCP does not always hold, we obtain weaker
results.
\begin{theorem}\label{thmmin_null}
Suppose that $1<p<\infty$, and let $Q_V$ be nonnegative on
$\core$. Let $\Omega_1\Subset \Omega$ be an open set, and let
$u\in C(\Omega\setminus \Omega_1)$ be a positive solution of the
equation $Q_V'(u)=0$ in $\Omega\setminus \overline{\Omega_1}$
satisfying $|\nabla u|\neq 0$ in $\Omega\setminus
\overline{\Omega_1}$.
\par
Then $u\in\M_{\Omega,\overline{\Omega_1}}$ if for every smooth
open set $\Omega_2$ satisfying $\Omega_1 \Subset \Omega_2\Subset
\Omega$, and an open set $B\Subset(\Omega\setminus
\overline{\Omega_2})$ there exists a sequence
$\{u_k\}\subset\core$,  $u_k\ge 0$, such that for all $k\in\N$,
$\int_B |u_k|^p\dx=1$, and \be
\lim_{k\to\infty}\int_{\Omega\setminus \overline{\Omega}_2}
L(u_k,u)\dx =0,\ee where $L$ is the Lagrangian given by
\eqref{piconeLag}.
\end{theorem}
\begin{proof}
Let $u$ be a positive solution of the equation $Q_V'(u)=0$ in
$\Omega\setminus \overline{\Omega}_1$ satisfying the theorem's
assumptions. Let $\Omega_2$ be an open set with smooth boundary
such that $\Omega_1 \Subset  \Omega_2\Subset \Omega$, and let $B$
an open set so that  $B\Subset(\Omega\setminus
\overline{\Omega_2})$.
\par
Suppose that $\{u_k\}\subset\core$ is a sequence of nonnegative
functions such that $\int_B |u_k|^p\dx=1$ for all $k\in\N$, and
\be\label{eq4.6} \lim_{k\to\infty}\int_{\Omega\setminus
\overline{\Omega_{2}}} L(u_k,u)\dx =0.\ee
As in the proof of \cite[Lemma~3.2]{ky3}, it follows  that $u_k\to
cu$ in $W^{1,p}_{\mathrm{loc}}
(\Omega\setminus\overline{\Omega_{1.5}})$, where $c>0$.

Let $E:W^{2,p}_{\mathrm{loc}}(\Omega \setminus \Omega_{2}) \to
W^{2,p}_{\mathrm{loc}}(\Omega)$ be the bounded extension operator,
constructed in the proof of Theorem~\ref{thmmin_null:2} (referring
here, as in \cite{E}, to the general case $1<p<\infty$ ):
$$
(Ew)(x):=w_\infty(x)\chi_\infty(x)+\sum_{i=0}^mw_i(x)\chi_i(x).
$$
Let $\tilde{u}=Eu$, $\tilde{u}_k=Eu_k$ and note that, since $u>0$,
one can always choose the neighborhoods $U_1,\dots,U_m$
sufficiently small so that $\tilde{u}>0$.

Let $f:=Q_V'(\tilde{u})$, and define $W:=f\tilde{u}^{1-p}$.
Clearly, $W$ has a compact support.  Since $|\nabla u|\neq 0$ in
$\Omega\setminus \overline{\Omega_1}$, a standard elliptic
regularity argument implies that $W\in
L^\infty_{\mathrm{loc}}(\Omega)$. It follows that $\tilde{u}$ is a
positive solution of the equation $Q_{V-W}'(u)=0$ in $\Omega$.
Moreover, by the continuity of $E$ and the continuity of $Q_{V-W}$
due \eqref{p<2}, it follows from \eqref{eq4.6} that \be
\lim_{k\to\infty} Q_{V-W}(\tilde{u}_k)=\lim_{k\to\infty}
Q_{V-W}(|\tilde{u}_k|)= \lim_{k\to\infty}\int_{\Omega}
L(|\tilde{u}_k|,\tilde{u})\dx =0.\ee Therefore,
\cite[Theorem~1.6]{ky3} implies that $|\tilde{u}_k|\to
c\tilde{u}$, and $\tilde{u}$ is a ground state of the functional
$Q_{V-W}$ in $\Omega$. Consequently, Theorem~\ref{thmmingr13}
implies that $\tilde{u}\in\M_{\Omega,\emptyset}$, and therefore
$u\in \M_{\Omega,\overline{\Omega_2}}\,$. Since $\Omega_2$ is an
arbitrary smooth open set satisfying  $\Omega_1 \Subset
\Omega_2\Subset \Omega$, we have $u\in
\M_{\Omega,\overline{\Omega_1}}\,$.
\end{proof}
\begin{remark}\label{remNecessity}{\em
Suppose that for all $V\in L^\infty_{\mathrm{loc}}(\Omega)$, any
positive solution of the equation $Q'_V(u)=0$ in $\Omega$
satisfying $u\in \M_{\Omega,\overline{\Omega_1}}$ for some smooth
open set $\Omega_1 \Subset \Omega$  is a global minimal solution
(this assumption seems to depend on the SCP, cf.
Proposition~\ref{prop_min}). Then the condition of
Theorem~\ref{thmmin_null} is also necessary.

Indeed, suppose that  $u\in C(\Omega\setminus\Omega_1)\cap
\M_{\Omega,\overline{\Omega_1}}$ satisfying $|\nabla u|\neq 0$ in
$\Omega\setminus \overline{\Omega_1}$. Let $\Omega_2$ be any open
set with smooth boundary such that $\Omega_1 \Subset
\Omega_2\Subset \Omega$.
Let $\tilde{u}$ be a positive smooth function in $\Omega$ such
that $\tilde{u}|_{\Omega\setminus\Omega_2}=u$. Let
$f:=Q_V'(\tilde{u})$, and define $W:=f\tilde{u}^{1-p}$. Recall
that $W$ has a compact support in $\Omega$, and by our assumption
$\tilde{u}$ is a global minimal solution of the equation
$Q'_{V-W}(u)=0$ in $\Omega$. So, by Theorem~\ref{thmmingr1},
$\tilde{u}$ is a ground state of the functional $Q_{V-W}$. Let
$\{u_k\}\subset \core$ be a null sequence for $Q_{V-W}$. So,
$\lim_{k\to\infty}Q_{V-W}(u_k)=0$, and for an open set
$B\Subset(\Omega\setminus \overline{\Omega_2})$,  we have  $\int_B
|u_k|^p\dx=1$ for all $k\in \mathbb{N}$. By Picone identity,
$$\lim_{k\to\infty}\int_\Omega
L(u_k,\tilde{u})\dx=\lim_{k\to\infty}Q_{V-W}(u_k)=0.$$  Since
$L(u_k,\tilde{u})\geq 0$ in $\Omega$, it follows that \be
\lim_{k\to\infty}\int_{\Omega\setminus \overline{\Omega}_2}
L(u_k,u)\dx =0.\ee\qed
 }\end{remark}
Finally, we formulate a sub-supersolution comparison principle for
our singular elliptic equation.
\begin{theorem}[Comparison Principle]\label{thm_comparison}
Suppose that $1<p<\infty$, and let $Q_V$ be nonnegative on
$\core$. Fix smooth open sets $\Omega_1\Subset \Omega_2\Subset
\Omega$. Let $u,v\in W^{1,p}_{\mathrm{loc}}(\Omega\setminus
\Omega_1)\cap C(\Omega\setminus \Omega_1)$ be, respectively, a
positive subsolution and a supersolution of the equation
$Q'_V(w)=0$ in $\Omega\setminus \overline{\Omega_1}$ such that
$u\leq v$ on $\partial \Omega_2$.
\par
Assume further that $Q'_V(u)\in
L^\infty_{\mathrm{loc}}(\Omega\setminus \Omega_1)$, $|\nabla
u|\neq 0$ in $\Omega\setminus \overline{\Omega_1}$, and that there
exist an open set $B\Subset(\Omega\setminus \overline{\Omega_2})$
and a sequence $\{u_k\}\subset\core$, $u_k\ge 0$, such that  \be
\int_B|u_k|^p\dx=1 \quad \forall k\geq 1, \mbox{ and }\quad
\lim_{k\to\infty}\int_{\Omega\setminus \overline{\Omega}_1}
L(u_k,u)\dx =0.\ee Then $u\leq v$ on $\Omega\setminus \Omega_2$.
\end{theorem}
\begin{proof} As in the proof of Theorem~\ref{thmmin_null}, let
$0<\tilde{u}\in W^{1,p}_{\mathrm{loc}}(\Omega)$ be an extension of
$u$ such that $\tilde{u}|_{\Omega\setminus\Omega_2}=u$. Let
$f:=Q_V'(\tilde{u})$, and define $W:=f\tilde{u}^{1-p}$. Clearly,
$W\leq 0$ in $\Omega\setminus \Omega_2$. By the proof of
Theorem~\ref{thmmin_null},  $\tilde{u}$ is a positive solution of
the equation $Q'_{V-W}(u)=0$ in $\Omega$ which is a global minimal
solution. On the other hand, $v$ is a positive supersolution of
the equation $Q'_{V-W}(w)=0$ in $\Omega\setminus \Omega_2$,
therefore $u\leq v$ on $\Omega\setminus \Omega_2$.
\end{proof}
\begin{remark} {\em Following Remark~\ref{remc1}, the normalization condition
$$\int_B|u_k|^p\dx=1$$
in Theorems~\ref{thmmin_null:2}, \ref{thm_comparison:2}, \ref{thmmin_null}
and \ref{thm_comparison} can be replaced by the condition
$$\int_B u_k\dx=1.$$}
\end{remark}
\begin{center}
{\bf Acknowledgments} \end{center}  The authors wish to thank
V.~Liskevich and P.~Tak\'{a}\v{c} for valuable discussions. The
work of Y.~P. was partially supported by the RTN network
``Nonlinear Partial Differential Equations Describing Front
Propagation and Other Singular Phenomena", HPRN-CT-2002-00274, and
the Fund for the Promotion of Research at the Technion.


\begin{thebibliography}{00} {\frenchspacing
%
%
\bibitem{Agmon84} S.~Agmon, Bounds on exponential decay of eigenfunctions of Schr\"odinger
operators, {\em in} ``Schr\"odinger Operators" (Como, 1984),
pp.~1--38, Lecture Notes in Math. 1159, Springer, Berlin, 1985.
%
\bibitem{AS} G.~Alessandrini, and M.~Sigalotti, Geometric properties of solutions
to the anisotropic $p$-Laplace equation in dimension two, {\em
Ann. Acad. Sci. Fenn. Math.} {\bf 26} (2001),  249--266.

\bibitem{AH1} W.~Allegretto, and Y.~X.~Huang, A Picone's identity for the
$p$-Laplacian and applications, {\em Nonlinear Anal.} {\bf 32}
(1998), 819--830.
%
\bibitem{AH2} W.~Allegretto, and Y.~X.~Huang, Principal eigenvalues and Sturm
comparison via Picone's identity, {\em  J. Differential Equations}
{\bf 156} (1999), 427--438.
%
\bibitem{CT} M.~Cuesta, and P.~Tak\'{a}\v{c}, A strong comparison principle for
positive solutions of degenerate elliptic equations, {\em
Differential Integral Equations} {\bf 13} (2000), 721--746.
%
\bibitem{DaS} L.~Damascelli, and B.~Sciunzi, Harnack inequalities,
maximum and comparison principles, and regularity of positive
solutions of $m$-Laplace equations, {\em Calc. Var. Partial
Differential Equations} {\bf 25} (2006), 139--159.
%
\bibitem{DS} J.~I.~Diaz, and J.~E.~Sa\'{a},  Existence et unicit\'{e} de
solutions positives pour certaines \'{e}quations elliptiques
quasilin\'{e}aires, {\em C. R. Acad. Sci. Paris Ser. I Math.} {\bf
305} (1987), 521--524.
%
\bibitem{DKN} P.~Dr\'{a}bek, A.~Kufner, and F.~Nicolosi, ``Quasilinear
Elliptic Equations with Degenerations and Singularities",  de
Gruyter Series in Nonlinear Analysis and Applications {\bf 5},
Walter de Gruyter \& Co., Berlin, 1997.
%
\bibitem{E} L.~C.~Evans, ``Partial Differential Equations",
Graduate Studies in Mathematics {\bf 19}, American Mathematical
Society, Providence, RI, 1998.
%
\bibitem{GS} J.~Garc\'{\i}a-Meli\'{a}n, and J.~Sabina de Lis, Maximum and comparison
principles for operators involving the $p$-Laplacian, {\em J.
Math. Anal. Appl.} {\bf 218} (1998), 49--65.
%
\bibitem{HKM} J.~Heinonen, T.~Kilpel\"{a}inen, and O.~Martio,
``Nonlinear Potential Theory of Degenerate Elliptic Equations",
Oxford Mathematical Monographs, Oxford University Press, New York,
1993.
%
\bibitem{LLM} V.~Liskevich, S.~Lyakhova, and  V.~Moroz, Positive
solutions to nonlinear $p$-Laplace equations with Hardy potential
in exterior domains, {\em J. Differential Equations}  {\bf 232}
(2007), 212--252.
%
\bibitem{LP} M.~Lucia, and S.~Prashanth, Strong comparison principle for solutions
of quasilinear equations, {\em Proc. Amer. Math. Soc.} {\bf 132}
(2004),  1005--1011.
%
\bibitem{M} J.~J.~Manfredi, Isolated singularities of $p$-harmonic
functions in the plane, {\em SIAM J. Math. Anal.} {\bf 22} (1991),
424--439.
%
\bibitem{MMP} M.~Marcus, V.~J.~Mizel, and Y.~Pinchover,  On the best constant
for Hardy's inequality in $\mathbb{R}^n$, {\em Trans. Amer. Math.
Soc.} {\bf 350} (1998), 3237--3255.
%
%
\bibitem{MP} \`{E}.~Mitidieri, and S.~I.~Pokhozhaev, Some generalizations of
Bernstein's theorem, {\em Differ. Uravn.} {\bf 38} (2002),
373--378; translation in {\em Differ. Equ.} {\bf 38} (2002),
392--397.
%
\bibitem{M86} M.~Murata, Structure of positive solutions to $(-\Delta
+V)u=0$ in $\mathbb{R}^n$, {\em Duke Math. J.} {\bf 53} (1986),
869--943.
%
\bibitem{PLiouv}  Y.~Pinchover, A Liouville-type theorem for Schr\"odinger
operators, {\em Comm. Math. Phys.} {\bf 272} (2007), 75--84.
%
\bibitem{P} Y.~Pinchover, Topics in the theory of positive solutions of
second-order elliptic and parabolic partial differential
equations, to appear in ``Spectral Theory and Mathematical
Physics: A Festschrift in Honor of Barry Simon's 60th Birthday",
eds. F~Gesztesy, et al., Proceedings of Symposia in Pure
Mathematics {\bf 76}, American Mathematical Society, Providence,
RI, 2007, 329--356.
%
\bibitem{aky} Y.~Pinchover, A.~Tertikas, and K.~Tintarev,
A Liouville-type theorem for the $p$-Laplacian with potential
term, to appear in {\em Ann.  Inst.  H. Poincar\'{e} (C) Anal. Non
Lin\'{e}aire} (2007), 12pp. doi:10.1016/j.anihpc.2006.12.004
%
%
\bibitem{ky2} Y.~Pinchover, and K.~Tintarev, Ground state alternative for singular
Schr\"{o}dinger operators, {\em J. Functional Analysis}, {\bf 230}
(2006), 65--77.
%
\bibitem{ky3} Y.~Pinchover, and K.~Tintarev, Ground state
alternative for $p$-Laplacian with potential term, {\em Calc. Var.
Partial Differential Equations} {\bf 28} (2007), 179--201.
%
\bibitem{PS} A.~Poliakovsky, and I.~Shafrir, Uniqueness of positive
solutions for singular problems involving the $p$-Laplacian, {\em
Proc. Amer. Math. Soc.} {\bf 133} (2005), 2549--2557.
%

\bibitem{Serrin1} J.~Serrin, Local behavior of solutions of quasi-linear
equations, {\em Acta Math.} {\bf 111} (1964), 247--302.
%
\bibitem{Serrin2} J.~Serrin, Isolated singularities of solutions of quasi-linear
equations, {\em Acta Math.} {\bf 113} (1965), 219--240.
%
\bibitem{TakTin} P.~Tak\'{a}\v{c}, and K.~Tintarev, Generalized minimizer solutions
for equations with the $p$-Laplacian and a potential term, to
appear in {\em  Proc. London Math. Soc.}.
%
\bibitem{T} P.~Tolksdorf,  Regularity for a more general class of
quasilinear elliptic equations, {\em J. Differential Equations}
{\bf  51} (1984), 126--150.
%
\bibitem{T1} P.~Tolksdorf, On the Dirichlet problem for quasilinear equations in
domains with conical boundary points, {\em Comm. Partial
Differential Equations} {\bf 8} (1983), 773--817.
%
\bibitem{Tr1} M.~Troyanov, Parabolicity of manifolds, {\em Siberian Adv. Math.} {\bf
9} (1999),  125--150.
%
\bibitem{Tr2} M.~Troyanov, Solving the $p$-Laplacian on manifolds, {\em Proc. Amer.
Math. Soc.} {\bf 128} (2000), 541--545.
%
\bibitem{V} L.~V\'{e}ron,  ``Singularities of Solutions of Second
Order Quasilinear Equations", Pitman Research Notes in Mathematics
Series, {\bf 353}. Longman, Harlow, 1996.
 }\end{thebibliography}
\end{document}